\documentclass[a4,10pt,twocolumn]{scrartcl}
\usepackage[top=2.5cm, bottom=2.5cm, left=1.6cm, right=1.6cm]{geometry}

\pdfoutput=1

\input{myStandardPreamble_arXiv.tex}
\usepackage[english]{babel}
\usepackage[utf8x]{inputenc}
\usepackage{enumitem}					

\usepackage{libertine}
\usepackage[T1]{fontenc}
\usepackage{mathpazo}

\usepackage[hang]{footmisc}
\setfnsymbol{wiley}

\makeatletter
\renewcommand{\@biblabel}[1]{[#1]\hfill}
\makeatother

\makeatletter
\let\NAT@parse\undefined
\makeatother

\usepackage[format=plain,			
            labelsep=period,		
            font=small,			
            labelfont=bf,			
            skip=5pt			
            ]{caption}
            
\usepackage[final,kerning=false,protrusion=false]{microtype}

\let\tmp\newinsert
\let\newinsert\newbox
\usepackage{floatrow}
\let\newinsert\tmp
\newfloatcommand{capbtabbox}{table}[][\FBwidth]

\usepackage{authblk}

\newtheorem{theorem}{Theorem}
\newtheorem{lemma}{Lemma}
\newtheorem{definition}{Definition}
\newtheorem{remark}{Remark}
\newtheorem{example}{Example}
\newtheorem{assumption}{Assumption}
\newtheorem{corollary}{Corollary}
\newtheorem{proposition}{Proposition}

\let\oldTheorem\theorem
\renewcommand{\theorem}{\oldTheorem\normalfont}
\let\oldLemma\lemma
\renewcommand{\lemma}{\oldLemma\normalfont}
\let\oldCorollary\corollary
\renewcommand{\corollary}{\oldCorollary\normalfont}
\let\oldDefinition\definition
\renewcommand{\definition}{\oldDefinition\normalfont}
\let\oldRemark\remark
\renewcommand{\remark}{\oldRemark\normalfont}
\let\oldAssumption\assumption
\renewcommand{\assumption}{\oldAssumption\normalfont}
\let\oldExample\example
\renewcommand{\example}{\oldExample\normalfont}

\addtokomafont{paragraph}{\itshape}

\colorlet{fictitiousColor}{green!45!black}    

\newcommand{\unitVec}[2][{}]{{}_{#1} e_{#2}}                		
\newcommand{\map}[3]{#1:#2 \rightarrow #3}						
\newcommand{\real}{\mathbb{R}}								
\newcommand{\pder}[2]{\frac{\partial #1}{\partial #2}}			
\newcommand{\oprocendsymbol}{\hbox{$\bullet$}} 					
\newcommand{\oprocend}{\relax\ifmmode\else\unskip\hfill\fi\oprocendsymbol}

\newcommand{\nIneq}{n_{\textup{ineq}}}                            	
\newcommand{\nEq}{n_{\textup{eq}}}	                        		
\newcommand{\setEq}{\mathcal{I}_{\textup{eq}}}					
\newcommand{\setIneq}{\mathcal{I}_{\textup{ineq}}}				
\newcommand{\dualIneq}{\lambda}                            		
\newcommand{\dualEq}{\nu}                                   		

\newcommand{\laplacian}{G}			                        		
\newcommand{\elLaplacian}{\MakeLowercase{\laplacian}}			
\newcommand{\G}{\mathcal{G}}								
\newcommand{\Adj}{\mathbf{A}}								
\newcommand{\AdjEl}{\mathbf{\MakeLowercase{\Adj}}}				
\newcommand{\degreeMat}{\mathrm{D}}							
\newcommand{\degreeMatEl}{\mathrm{\MakeLowercase{\degreeMat}}}		
\newcommand{\pathLeft}{\langle}								
\newcommand{\pathRight}{\rangle}							
\newcommand{\pathSep}{|}									

\newcommand{\state}{z}									
\newcommand{\dimState}{N}									
\newcommand{\seqParam}{\sigma}								

\newcounter{MyAssumptionCounter}
\renewcommand{\emph}{\textit}

\newcommand{\simon}[1]{{\color{orange} #1}}

\newcommand{\margin}[1]{\marginpar{\color{red}\tiny\ttfamily#1}}

\renewcommand{\margin}[1]{\marginpar{\color{red}\tiny\ttfamily}}

\renewcommand\simon[1]{#1}

\newboolean{longVersion}
\setboolean{longVersion}{true}

\title{\LARGE \textbf
On the Lie bracket approximation approach to distributed optimization: Extensions and limitations
\footnote{This article is an extended version of \cite{mic2018extensions} additionally including proofs for all results, a discussion of the choice of vector fields after \cref{remarkStructureBoundedVecFields} and the overview in \cref{tableVectorFields}.}
}

\setkomafont{section}{\Large}
\setkomafont{subsection}{\large}
\setkomafont{subsubsection}{\itshape}

\begin{document}

\date{}
\author[1]{Simon Michalowsky}
\author[2]{Bahman Gharesifard}
\author[1]{Christian Ebenbauer}
\affil[1]{Institute for Systems Theory and Automatic Control, University of Stuttgart, Germany \protect\\ \texttt{\small $\lbrace$michalowsky,ce$\rbrace$@ist.uni-stuttgart.de}}
\affil[2]{Department of Mathematics and Statistics, Queen's University, Canada \protect\\ \texttt{\small bahman@mast.queensu.ca} \protect\\[1em]}

\maketitle

\begin{abstract}
\textbf{Abstract.} We consider the problem of solving a smooth convex optimization problem with 
	equality and inequality constraints in a distributed fashion. Assuming that 
	we have a group of agents available capable of communicating over a communication network
	described by a time-invariant directed graph, we derive distributed continuous-time agent dynamics
	that ensure convergence to a neighborhood of the optimal solution of the optimization problem. 
	Following the ideas introduced in our previous work, 
	we combine saddle-point dynamics with Lie bracket approximation techniques. 
	While the methodology was previously limited to linear constraints and objective functions
	given by a sum of strictly convex separable functions, we extend these
	result here and show that it applies to a 
	very general class of optimization problems under mild assumptions on the 
	communication topology.
\end{abstract}

\section{INTRODUCTION}\label{secIntroduction}
Over the last decades, distributed optimization has been an active area of research 
with high practical relevance, see, e.g., \cite{boyd2011distributed, bullo2009robotic, zhao2015microgrids} for applications. 
In these type of problems, the goal is to cooperatively solve an optimization problem
using a group of agents communicating over a network. While discrete-time
algorithms for distributed optimization constitute the majority in the existing
literature, we focus on continuous-time algorithms which have regained interest 
in the last decades, 
\cite{feijer2010network, wang2011control, durr2013saddle, gharesifard2014distributed, niederlaender2015distributed, touri2016saddle}. 
These algorithms require strong assumptions either on the structure of the 
optimization problem or on the communication network.
Recently, a novel approach to continuous-time distributed optimization based
on Lie bracket approximations has been proposed in~\cite{ebenbauer2017directed},
\cite{mic2017extremum} that has the potential to relax these assumptions.

In the present work we want to extend these results and show that the approach
can be applied to a large class of optimization problems under mild assumptions
on the communication network. While in~\cite{ebenbauer2017directed}, \cite{mic2017extremum}, \cite{mic2017opt}
only optimization problems with {linear constraints} and objective
functions in the form of a sum of separable functions were considered,
in the present work we enhance the methodology to general convex optimization problems.
The main idea of the approach is to use Lie bracket approximation techniques
to find distributed approximations of non-distributed saddle-point dynamics.
While in the previous works only certain Lie bracket approximations were used,
we further show here that a whole class is applicable.

\section{PRELIMINARIES}\label{secPreliminaries}
\subsection{Notation}
We denote by $ \mathbb{R}^n $ the set of $n$-dimensional real
vectors and further write $ \mathcal{C}^p $, $ p \in \mathbb{N} $, for 
the set of $p$-times continuously differentiable real-valued functions.
The gradient of a function $ f: \real^n \to \real $, $ f \in \mathcal{C}^1 $, 
with respect to its argument $x\in \real^n$,
will be denoted by $  \map{\nabla_x f}{\real^n}{\real^n} $;
we often omit the subscript, if it is clear from the context.
We denote the $ (i,j) $th entry of a matrix $ A \in \mathbb{R}^{n \times m} $
by $ a_{ij} $, and sometimes denote $ A $ by $ A = [a_{ij}] $. We use $ \unitVec{i} $ to denote the vector 
with the $ i $th entry equal to $ 1 $ and all other entries equal to~$0$.
We do not specify the dimension of $ \unitVec{i} $ but expect it to be
clear from the context.
For a vector $ \lambda \in \mathbb{R}^n $ we let $ \textup{diag}(\lambda) \in \mathbb{R}^{n \times n} $
denote the diagonal matrix whose diagonal entries are the entries of $ \lambda $.
Given two continuously differentiable vector fields 
$ \map{f_1}{\real^n}{\real^n} $ and $ \map{f_2}{\real^n}{\real^n} $,
the Lie bracket of $ f_1 $ and $ f_2 $ evaluated at $x$ is defined to be
\begin{align}
	[ f_1, f_2 ](x) := \pder{f_2}{x}(x) f_1(x) - \pder{f_1}{x}(x)  f_2(x). \label{eqDefLieBracket}
\end{align}
With a slight abuse of notation we sometimes also write $ [ f_1(x), f_2(x) ] $.
For a vector $ x = [ x_1, x_2, \dots, x_n ]^\top\in\mathbb{R}^n $ and
a finite set $ \mathcal{S} \subseteq \lbrace 1,2,\dots,n \rbrace $, we denote
by $ x_{\mathcal{S}} \in \mathbb{R}^{\vert \mathcal{S} \vert} $ the ordered stacked vector of all $ x_i $ with
$ i \in \mathcal{S} $. For example, if $ n = 5 $ and $ \mathcal{S} = \lbrace 2,4 \rbrace $,
then $ x_{\mathcal{S}} = [ x_2, x_4 ]^\top $.

\subsection{Basics on graph theory}
We recall some basic notions on graph theory, and refer the reader 
to~\cite{biggs1993algebraic} or other standard references for more information.
A directed graph (or simply digraph) is an ordered pair $ \G = ( \mathcal{V}, \mathcal{E} ) $,
where $ \mathcal{V} = \lbrace v_1, v_2, \dots, v_n \rbrace $
is the set of nodes and $ \mathcal{E} \subseteq \mathcal{V} \times \mathcal{V} $
is the set of edges, i.e. $ (v_i,v_j) \in \mathcal{E} $ if there 
is an edge from node $ v_i $ to $ v_j $. In our setup the
edges encode to which other agents some agent has access to, i.e. 
$ (v_i, v_j) \in \mathcal{E} $ means that node $ v_i $
receives information from node $ v_j $. We say that node $ v_j $ is an out-neighbor of
node $v_i$ if there is an edge from node $ v_i $ to node $ v_j $. 
\ifthenelse{\boolean{longVersion}}{
The adjacency matrix $ \Adj = [ \AdjEl_{ij} ]  \in \real^{n\times n} $ associated to $ \mathcal{G} $ is defined as
\begin{align}
	\AdjEl_{ij}
	&=
	\begin{cases}
		1 & \text{if } i \neq j \text{ and } (v_i,v_j) \in \mathcal{E},  \\
		0 & \text{otherwise}.
	\end{cases} \label{eqPreliminariesDefAdjacencyMat}
\end{align}
We also define the out-degree matrix  $ \degreeMat = [ \degreeMatEl_{ij} ] $ associated to $ \G $ as
\begin{align}
	\degreeMatEl_{ij}
	&=
	\begin{cases}
		\sum_{k=1}^{n} \AdjEl_{ik} & \text{if } i = j \\
		0 \hspace*{1.3cm}                   & \text{otherwise.}
	\end{cases}
\end{align}
}{
The adjacency matrix $ \Adj = [ \AdjEl_{ij} ]  \in \real^{n\times n} $ associated to $ \mathcal{G} $ is defined as
\begin{align}
	\AdjEl_{ij}
	&=
	\begin{cases}
		1 & \text{if } i \neq j \text{ and } (v_i,v_j) \in \mathcal{E},  \\
		0 & \text{otherwise}.
	\end{cases} \label{eqPreliminariesDefAdjacencyMat}
\end{align}
We also define the out-degree matrix  $ \degreeMat = [ \degreeMatEl_{ij} ] $ associated to $ \G $ as
\begin{align}
	\degreeMatEl_{ij}
	&=
	\begin{cases}
		\sum_{k=1}^{n} \AdjEl_{ik} & \text{if } i = j \\
		0 \hspace*{1.3cm}                   & \text{otherwise.}
	\end{cases}
\end{align}
}
Finally, we call $ \laplacian = \degreeMat - \Adj = [ \elLaplacian_{ij} ] \in \mathbb{R}^{n \times n} $ the Laplacian of $ \mathcal{G} $. 
A directed path in $ \mathcal{G} $ is a sequence of nodes connected by edges
and we write $ p_{i_1 i_r} = \pathLeft v_{i_1} \pathSep v_{i_2} \pathSep \dots \pathSep v_{i_r} \pathRight $
for a path from node $ v_{i_1} $ to node $ v_{i_r} $. We say that a path
is simple if the sequence contains no node more than once.

\section{PROBLEM SETUP}\label{secProblemSetup}
Consider the following convex optimization problem
\begin{align}
	\begin{split}
	\min        \quad & F(x) \\
	\text{s.t.} \quad & a_i(x) = 0, \quad i \in \setEq \subseteq \lbrace 1,2,\dots,n \rbrace, \\
	                  & c_i(x) \leq 0, \quad i \in \setIneq \subseteq \lbrace 1,2,\dots,n \rbrace,
	\end{split}
	\label{eqOptimizationProblem}
\end{align}
where $ F: \mathbb{R}^n  \to \mathbb{R} $, $ a_i: \mathbb{R}^n \to \mathbb{R}^{{\nEq}_i} $, 
$ c_i: \mathbb{R}^n \to \mathbb{R}^{{\nIneq}_i} $, $ {\nEq}_i, {\nIneq}_i \geq 1 $,
$ F \in \mathcal{C}^2 $ is strictly convex, the functions $ a_i,i \in \setEq, $
are affine and the $ c_i \in \mathcal{C}^2,i\in\setIneq $, are convex.  
We assume further that $ F,a_i,c_i$ are such that the feasible set of~\eqref{eqOptimizationProblem} 
is non-empty and that the problem has a unique solution.

Our goal is to design continuous-time optimization algorithms that
converge to an arbitrarily small neighborhood of the unique global
optimizer of~\eqref{eqOptimizationProblem} and that can be implemented
in a distributed fashion. More precisely, we assume that we have
a group of $n$ agents available, each capable of interchanging information
over a communication network {described} by a directed graph $ \mathcal{G} = ( \mathcal{V}, \mathcal{E} ) $
with graph Laplacian $ \laplacian = [ \elLaplacian_{ij} ] $,
where $ \mathcal{V} = \lbrace v_1, v_2, \dots, v_n \rbrace $ is a set of
$n$ nodes and $ \mathcal{E} \subseteq \mathcal{V} \times  \mathcal{V} $ 
is the edge set between the nodes. In the present setup, each node~$v_i$
represents an agent and the edges define the existing communication links
between the agents, i.e., if there is an edge from node~$i$ to node~$j$
then agent~$i$ has access to the information provided by agent~$j$. 
We then say that an algorithm is distributed if each agent only uses its
own information as well as that provided by its out-neighboring agents.

Let $ L: \mathbb{R}^n \times \mathbb{R}^{\nEq} \times \mathbb{R}^{\nIneq} \to \mathbb{R} $,
$ \nEq = \sum_{i \in \setEq} {\nEq}_i $, $ \nIneq=~\sum_{i \in \setIneq} {\nIneq}_i $, 
denote the Lagrangian associated to~\eqref{eqOptimizationProblem}, i.e.,
\begin{align}
	L(x,\nu,\lambda) = F(x) + \dualEq^\top a(x) + \dualIneq^\top c(x),
\end{align}
where $ a = [ a_i ]_{i \in \setEq}, c = [ c_i ]_{i \in \setIneq} $ are the
stacked vectors of all $ a_i $ and $ c_i $, respectively, and $ \dualEq \in \mathbb{R}^{\nEq} $,
$ \dualIneq \in \mathbb{R}^{\nIneq} $ are the associated Lagrange
multipliers. {In the sequel,} we assume that the state of the $i$th agent 
comprises of $x_i$ as well as the dual variables associated to the
constraints $a_i$, $c_i$, given they exist. 
Without loss of generality we then put the following assumption on the indexing
of the constraints:
\begin{assumption}
For any $ i \in \setEq $ and any $ j \in \setIneq ${,} there exists an $ x~=[ x_1, x_2, \dots, x_n ]^\top \in\mathbb{R}^n $ such that
	$ \tfrac{\partial a_i}{\partial x_i}(x) \neq 0 $ and 
	$ \tfrac{\partial c_j}{\partial x_j}(x) \neq 0 $.
\end{assumption}
In a nutshell, this assumption New{guarantees} that the $i$th constraints $a_i$, $c_i$
are both functions of $x_i$. 

It is well-known that if the Lagrangian has New{a saddle point}
$ (x^\star,\dualEq^\star,\dualIneq^\star) $, then $ x^\star $ is an optimizer
of~\eqref{eqOptimizationProblem}. We say that a point $ (x^\star,\dualEq^\star,\dualIneq^\star) $
is a saddle point if for all $ x \in \mathbb{R}^n $, $ \dualEq \in \mathbb{R}^{\nEq} $, {$ \dualIneq \in \mathbb{R}^{\nIneq}_+ $},
we have
\begin{align}
	L(x^\star,\dualEq,\dualIneq) \leq L(x^\star,\dualEq^\star,\dualIneq^\star) \leq L(x,\dualEq^\star,\dualIneq^\star).
\end{align}
The following saddle-point dynamics adapted from \cite{duerr2012saddlePoint} is known to converge
to a saddle point of the Lagrangian (see~\cite{mic2017opt} for a proof), thus providing a 
solution to~\eqref{eqOptimizationProblem}
\begin{subequations}
\ifthenelse{\boolean{longVersion}}{
\begin{align}
	\dot{x}         &= - \nabla_x L(x,\dualEq,\dualIneq)                                  \nonumber \\
	                &= -\nabla F(x) - \tfrac{\partial a}{\partial x}(x)^\top \dualEq -  \tfrac{\partial c}{\partial x}(x)^\top \dualIneq &   \label{eqSPAa} \\
	\dot{\dualEq}   &= \nabla_{\dualEq} L(x,\dualEq,\dualIneq)                            \nonumber \\
	                &= a(x)  & \label{eqSPAb} \\
	\dot{\dualIneq} &= \textup{diag}(\dualIneq) \nabla_{\dualIneq} L(x,\dualEq,\dualIneq) \nonumber \\
	                &= \textup{diag}(\dualIneq) c(x). &  \label{eqSPAc}
\end{align}}{
\begin{align}
	\dot{x}         &= -\nabla F(x) - \tfrac{\partial a}{\partial x}(x)^\top \dualEq -  \tfrac{\partial c}{\partial x}(x)^\top \dualIneq &   \label{eqSPAa} \\
	\dot{\dualEq}   &= a(x)  & \label{eqSPAb} \\
	\dot{\dualIneq} &= \textup{diag}(\dualIneq) c(x). &  \label{eqSPAc}
\end{align}}
\label{eqSPA}
\end{subequations}
However,~\eqref{eqSPA} is in general not distributed in the aforementioned
sense, since the right-hand side of~\eqref{eqSPA} is not composed only 
of \emph{admissible} vector fields, i.e., vector fields that can be computed 
locally by the nodes. For example, if $ \elLaplacian_{12} \neq 0 $ and $ \elLaplacian_{13} = 0 $, then the 
vector field $ [ \unitVec{1}^\top x_2, 0, 0 ]^\top $ is admissible for~\eqref{eqSPA}, while 
$ [ \unitVec{1}^\top x_3, 0, 0 ]^\top $ is not. 

Recently, a novel approach to distributed optimization has been proposed
that employs Lie bracket approximation techniques to derive distributed
approximations of~\eqref{eqSPA}. The idea is to write the right-hand side
of~\eqref{eqSPA} by means of Lie brackets of {admissible} vector fields.
If we have achieved to rewrite \eqref{eqSPA} in this form, i.e., we have
\begin{align}
	\dot{\state} = \sum\limits_{B \in \mathcal{B}} v_B B(\state), \label{eqProbSetupLieBracketSystem}
\end{align}
where $ z = [ x^\top, \dualEq^\top, \dualIneq^\top ]^\top $, $ v_B \in \mathbb{R} $ and
$ \mathcal{B} $ is a set of Lie brackets of admissible vector fields 
$ \lbrace \phi_1, \phi_2, \dots, \phi_M \rbrace $, $ \phi_i: \mathbb{R}^{n + \nEq + \nIneq} \to \mathbb{R}^{n + \nEq + \nIneq} $,
then we can employ Lie bracket approximation techniques from \cite{liu1997approximation}, \cite{sussmann1991limits} to derive
distributed approximations of \eqref{eqProbSetupLieBracketSystem}.
More precisely, we can find a family of functions $ u_i^{\seqParam} : \mathbb{R} \to \mathbb{R} $
parametrized by $ \seqParam > 0 $ such that the trajectories of
\begin{align}
	\dot{\state}^{\seqParam} = \sum\limits_{i=1}^{M} \phi_i(\state) u_i^{\seqParam}(t) \label{eqProbSetupApproximatingSystem}
\end{align}
uniformly converge to those of~\eqref{eqProbSetupLieBracketSystem}
as $ \seqParam $ increases, given~\eqref{eqProbSetupLieBracketSystem}
and~\eqref{eqProbSetupApproximatingSystem} are initialized equally. A {general} algorithm to compute suitable
functions $ u_i^{\sigma} $ is presented in~\cite{liu1997approximation}
{and we will present a modified version thereof tailored to the problem at hand
in \cite{mic2017opt}}.
In the present paper we do not discuss the second step of how 
to design these functions $ u_i^\seqParam $ but focus
on rewriting the right-hand side of~\eqref{eqSPA}
in terms of Lie brackets of admissible vector fields.

\section{MAIN RESULTS}\label{secMainResults}
Consider the saddle-point dynamics~\eqref{eqSPA}
and observe that the right-hand side is a sum of vector fields of the form
$ \unitVec{i} \psi(x) $ and $ \unitVec{i} \state_j \psi(x) $,
$ i,j = 1,2,\dots,n+\nEq+\nIneq $, $ \psi:~\mathbb{R}^n \to \mathbb{R} $, 
where 
\begin{align} 
	\state = [ x^\top, \dualEq^\top, \dualIneq^\top ]^\top \in \mathbb{R}^{\dimState}, \quad \dimState = n+\nEq+\nIneq,  \label{eqDefCompleteState} 
\end{align}
is the complete state and $ \unitVec{i} \in \mathbb{R}^{\dimState} $ is the $i$th
unit vector. These vector fields might either be admissible or
not, depending on the communication graph as well as the problem 
structure. In the following, we wish to discuss how to write vector 
fields of this form by means of Lie brackets of admissible vector fields.
\simon{For the purpose of notation}, for each $ i = 1,2,\dots,n $, we define the
index set 
\begin{flalign}
	&\bar{\mathcal{I}}(i) := \lbrace i \rbrace \cup \lbrace j = n + \sum\limits_{k=1}^{i-1} {\nEq}_k + \ell, \ell = 1,2,\dots,{\nEq}_{i}  \rbrace & \nonumber \\
	&\cup \lbrace j = n + \nEq + \sum\limits_{k=1}^{i-1} {\nIneq}_k + \ell, \ell = 1,2,\dots,{\nIneq}_{i}  \rbrace \hspace*{-2em} & \label{eqDefIndexSetOld}
\end{flalign}
associating the components of the complete state $ \state $ to the 
$i$th agent meaning that $ \state_j $ is part of the state of agent $i$
for all $ j \in \bar{\mathcal{I}}(i) $. Hence, the state vector of 
the $i$th agent is given by $ \state_{\bar{\mathcal{I}}(i)} $.
Based on this, for all $ j = 1,2,\dots,\dimState $, we define
\begin{align}
	\mathcal{I}(j) = \bar{\mathcal{I}}(i) \label{eqDefIndexSet}
\end{align}
for some $ i $ such that $ j \in \bar{\mathcal{I}}(i) $, i.e., $ \mathcal{I}(j) $
is the set of all indices which are associated to the same agent as
the $j$th index. Note that $ \bar{\mathcal{I}}(j) = \mathcal{I}(j) $, for all $ j = 1,2,\dots,n $.

\subsection{Lie brackets of admissible vector fields}
In the following, we first want to discuss which kind of vector fields
can be written in terms of Lie brackets of admissible vector fields.
For $ i,j = 1,2,\dots, \dimState $, define
\begin{align}
	h_{i,j}(\state) = \unitVec{j} f_j(\state_{\mathcal{I}(j)},\state_{\mathcal{I}(i)}), \label{eqDefVectorFields}
\end{align}
where $ f_j: \mathbb{R}^{\vert \mathcal{I}(j) \vert} \times \mathbb{R}^{\vert \mathcal{I}(i) \vert} \to \mathbb{R} $,
$ f_j \in \mathcal{C}^1 $. Observe that $ h_{i,j} $ is admissible
if and only if there exist $ \ell, k \in \lbrace 1,2,\dots,n \rbrace $ such that $ i \in \mathcal{I}(\ell) $,
$ j \in \mathcal{I}(k) $ and $ \elLaplacian_{k\ell} \neq 0 $. In the next New{result}
we consider Lie brackets of admissible vector fields of the form~\eqref{eqDefVectorFields}.
\begin{lemma}\label{lemmaGeneralFormula}
Consider a graph $ \mathcal{G} $ of $n$ nodes and let $ p_{i_1 i_r}~=~\pathLeft v_{i_1} \pathSep v_{i_2} \pathSep \dots \pathSep v_{i_r} \pathRight $
	denote a simple path in $ \mathcal{G} $ from $ v_{i_1} $ to $ v_{i_r} $. 
	Then, for any $ j_k \in \mathcal{I}(i_k) $, $ k = 1,2,\dots, r $, we have
	\begin{flalign}
		&\bigg[ h_{j_{r},j_{r-1}}, \Big[ h_{j_{r-1},j_{r-2}}, \big[ \dots, [ h_{j_{3},j_{2}}, h_{j_{2},j_{1}} ] \dots \big] \Big] \bigg](\state) \label{eqGeneralFormula}  & \\
		&=
		\unitVec{j_1} f_{j_{r-1}}( \state_{\mathcal{I}(j_{r-1})}, \state_{\mathcal{I}(j_r)} ) \prod\limits_{k=1}^{r-2} \tfrac{\partial f_{j_k}}{\partial \state_{j_{k+1}}}(\state_{\mathcal{I}(j_k)},\state_{\mathcal{I}(j_{k+1})}) 
		\nonumber 
	\end{flalign}
	and the left-hand side is a Lie bracket of admissible vector fields.
	\oprocend
\end{lemma}
A proof is given in \ifthenelse{\boolean{longVersion}}{\cref{appProofLemmaGeneralFormula}}{\cite{mic2018extensionsArxiv}}.
By the above Lemma, each non-admissible vector field that takes 
the same form as the right-hand side of~\eqref{eqGeneralFormula}
can be written in terms of a Lie bracket of admissible vector fields.
It is worth mentioning that this does not classify the whole set
of vector fields that can be written as a Lie bracket of admissible
vector fields since we limited ourselves to a single path.
We next discuss a special case that is of particular importance
for the application at hand.
\begin{proposition}\label{propFormulaSpecialCase}
Consider a graph $ \mathcal{G} $ of $n$ nodes and let $ p_{i_1 i_r} = \pathLeft v_{i_1} \pathSep v_{i_2} \pathSep \dots \pathSep v_{i_r} \pathRight $
	denote a simple path in $ \mathcal{G} $ from $ v_{i_1} $ to $ v_{i_r} $.
	Let $ j_k \in \mathcal{I}(i_k), k=1,2,\dots,r $, be any set of indices
	and suppose that
	\begin{align}
		f_{j_k}(\state_{\mathcal{I}(j_k)},\state_{\mathcal{I}(j_{k+1})}) = f_{j_k}^{(1)}(\state_{\mathcal{I}(j_k)}) f_{j_k}^{(2)}(\state_{\mathcal{I}(j_{k+1})})
		\label{eqPropFormulaSpecialCaseStructureFunctions}
	\end{align}
	with
	\begin{align}
		& f_{j_k}^{(1)}(\state_{\mathcal{I}(j_k)}) \tfrac{ \partial f_{j_{k-1}}^{(2)} }{ \partial \state_{j_k} } (\state_{\mathcal{I}(j_k)}) = 1 \label{eqPropConstraintFunctions}
	\end{align}
	for all {$ \state_{\mathcal{I}(j_k)} \in \mathbb{R}^{\vert \mathcal{I}(j_k) \vert} $}, $ k = 2,3,\dots,r-1 $. Then
	\begin{align}
		\begin{split}
		\bigg[ h_{j_{r},j_{r-1}}, \Big[ h_{j_{r-1},j_{r-2}}, \big[ \dots, [ h_{j_{3},j_{2}}, h_{j_{2},j_{1}} ] \dots \big] \Big] \bigg](\state) \\
		=
		\unitVec{j_1} f_{j_1}^{(1)}(\state_{\mathcal{I}(j_1)}) f_{j_{r-1}}^{(2)}( \state_{\mathcal{I}(j_r)} )
		\end{split}\label{eqPropFormulaSpecialCase}
	\end{align}
	for all $ \state \in \mathbb{R}^{\dimState} $
	and the left-hand side is a Lie bracket of admissible vector fields.
	\oprocend
\end{proposition}
A proof is given in \ifthenelse{\boolean{longVersion}}{\cref{appProofPropFormulaSpecialCase}}{\cite{mic2018extensionsArxiv}}.
{In view of~\eqref{eqGeneralFormula}, the constraint~\eqref{eqPropConstraintFunctions}
ensures that all terms depending on $ \state_{\mathcal{I}(j_k)} $, $ k = 2,3,\dots,r-1 $,
cancel out.}
Equation \eqref{eqPropFormulaSpecialCase} is of particular interest since
the non-admissible vector fields often take the form as its right-hand side. 
According to~\eqref{eqPropConstraintFunctions}, there exists a whole class of vector fields
$ h_{j_{k},j_{k-1}} $, or equivalently, functions $ f_{j_k} $, such 
that~\eqref{eqPropFormulaSpecialCase} holds. A particularly
simple choice that has been utilized in the previous works~\cite{ebenbauer2017directed},
\cite{mic2017extremum} is to take
\begin{subequations}
\begin{align}
	f_{j_k}^{(1)}(\state_{\mathcal{I}(j_k)}) &= 1,                    && k = 2,3,\dots,r-1, \\ 
	f_{j_k}^{(2)}(\state_{\mathcal{I}(j_{k+1})}) &= \state_{j_{k+1}}, && k = 1,2,\dots,r-2;
\end{align}\label{eqVectorFieldsSimpleChoice}
\end{subequations}
hence leading to $ h_{j_{k+1},j_{k}}(\state_{\mathcal{I}(j_k)},\state_{\mathcal{I}(j_{k+1})}) = \unitVec{j_k} \state_{j_{k+1}} $.
However, in view of~\eqref{eqProbSetupApproximatingSystem}
where the $ \phi_i $ are given by admissible vector fields of the 
form~\eqref{eqDefVectorFields}, it is often desired that the admissible vector
fields have certain properties such as boundedness {in order to
simplify the calculation of the approximating inputs $u_i^{\seqParam} $
and improve the transient behavior of~\eqref{eqProbSetupApproximatingSystem}}. However,
as we see next, it is not possible to render all vector fields
bounded.
\begin{lemma}
Suppose that all assumptions from \cref{propFormulaSpecialCase} are fulfilled.
	Then there exists no set of bounded vector fields {$ h_{j_k,j_{k-1}} \in \mathcal{C}^1 $}
	such that~\eqref{eqPropFormulaSpecialCase}
	holds.
	\oprocend
\end{lemma}
\ifthenelse{\boolean{longVersion}}{
\begin{proof}
Suppose there exists a set of bounded vector fields such 
	that~\eqref{eqPropFormulaSpecialCase} holds. Then 
	$ f_{j_k} $ is a bounded function for all $ k~=~1,2,\dots,r $.
	{By~\eqref{eqPropFormulaSpecialCaseStructureFunctions},
	$ f_{j_k}^{(1)} $ is also bounded for all $k~=~1,2,\dots,r$
	{since $\mathcal{I}(j_k)$ and $\mathcal{I}(j_{k+1})$ are disjunct}.}
	However, by \eqref{eqPropConstraintFunctions}, 
	\begin{align}
		\tfrac{ \partial f_{j_{k-1}}^{(2)} }{ \partial \state_{j_k} } (\state_{\mathcal{I}(j_k)})
		= \tfrac{1}{f_{j_k}^{(1)}(\state_{\mathcal{I}(j_k)})}
	\end{align}
	for $ k = 2,\dots,r-1 $; hence $ \tfrac{ \partial f_{j_{k-1}}^{(2)} }{ \partial \state_{j_k} } $
	is bounded away from zero, i.e., there exists some constant $ \beta > 0 $ such
	that for all $ \state_{\mathcal{I}(j_k)} \in \mathbb{R}^{\vert \mathcal{I}(j_k) \vert} $ we have 
	\begin{align}
		\vert \tfrac{ \partial f_{j_{k-1}}^{(2)} }{ \partial \state_{j_k} }(\state_{\mathcal{I}(j_k)}) \vert \geq \beta.
	\end{align}
	{Note that $ h_{j_k,j_{k-1}} \in \mathcal{C}^1 $ if and only if
	$ f_{j_{k-1}}^{(1)}, f_{j_{k-1}}^{(2)} \in \mathcal{C}^1 $.}	
	{Thus, $ f_{j_{k-1}}^{(2)} $ is strictly
	monotone} in $ \state_{j_k} $ which contradicts the boundedness assumption,
	thus concluding the proof.
\end{proof}
}{}
\begin{remark}\label{remarkStructureBoundedVecFields}
The same holds true if we do not assume the structure~\eqref{eqPropFormulaSpecialCaseStructureFunctions}. 
	In fact, the structure is required for all other variables except
	of $ \state_{\mathcal{I}(j_1)} $, $ \state_{\mathcal{I}(j_{r})} $
	to cancel out in~\eqref{eqPropFormulaSpecialCase}.
\oprocend
\end{remark}
\ifthenelse{\boolean{longVersion}}{As we see from the proof, each bounded vector field in~\eqref{eqPropFormulaSpecialCase}
leads to another unbounded vector field. Hence, at most half
of the vector fields $ h_{j_k,j_{k-1}} $ in~\eqref{eqPropFormulaSpecialCase}
can be bounded. In particular, we can choose the functions $ f_{j_k} $, 
$ k~=~1,2,\dots,r-1$, {as follows to guarantee} 
that~\eqref{eqPropConstraintFunctions} holds:
\begin{subequations}
\begin{flalign}
	f_{j_k}^{(1)}(x) 
	&=
	\begin{cases}
		\cos(\alpha\big(x-d)\big)            & \text{if } k \text{ is even}, \\
		\tfrac{1}{\cos(\alpha(x-d))} & \text{if } k \text{ is odd}, k \neq 1,
	\end{cases} \\
	f_{j_k}^{(2)}(x)
	&=
	\begin{cases}
		\sin\big(\alpha(x-d)\big)                            & \hspace*{-0.4em} \text{if } k \text{ is even}, k\neq r-1 , \\
		\smallint \tfrac{1}{\cos(\alpha(x-d)) } \mathrm{d} x & \hspace*{-0.4em} \text{if } k \text{ is odd}, k\neq r-1 ,
	\end{cases} \hspace*{-2em} &
\end{flalign}
\end{subequations}
where $ \alpha \neq 0 $, $ d \in \mathbb{R} $. However, this choice will lead
to functions $ f_{j_k} $ that are not globally continuous but only well-defined in 
the interval $ ( d-\tfrac{\pi}{2\alpha}, d+\tfrac{\pi}{2\alpha} ) $; hence
we need to choose $d$ appropriately and $ \alpha $ sufficiently small.
By that choice, all functions $ f_{j_k} $ with even $ k $ are bounded
while all functions $ f_{j_k} $ with odd $k$ are unbounded. }{}

\ifthenelse{\boolean{longVersion}}{}{The proof is given in \cite{mic2018extensionsArxiv}.}
We next discuss how we can make use of the previous results to rewrite
more general vector fields.
While \cref{lemmaGeneralFormula} enables us to write products of functions
of variables of nodes which lie on the same path as Lie brackets of admissible
vector fields, we cannot directly use this result to rewrite functions
which do not fulfill this property. 
To make this clearer, consider the following example:
\begin{example}\label{exampleGeneralFormula}
Consider the graph shown in \cref{figGraphExampleProduct} and assume for the 
	sake of simplicity that $ \mathcal{I}(j) = j $ for $ j = 1,\dots,5 $ and $ \state \in \mathbb{R}^5 $.
	By~\eqref{eqGeneralFormula} we can write non-admissible vector fields of the form
	$ \unitVec{1} \varphi_1(\state_2) \varphi_2(\state_3) $ as well
	as $ \unitVec{1} \varphi_3(\state_4) \varphi_4(\state_5) $ and sums thereof in 
	terms of Lie brackets of admissible vector fields as long as
	$ \varphi_1, \varphi_3 $ admit an analytic expression of their
	antiderivatives $ \smallint \varphi_1(\state_2) \mathrm{d}\state_2 $,
	$ \smallint \varphi_3(\state_4) \mathrm{d}\state_4 $. In fact, with 
	\begin{subequations}
	\begin{align}
		h_{2,1}(\state) &= \unitVec{1} \smallint \varphi_1(\state_2) \mathrm{d}\state_2, \quad h_{3,2}(\state) = \unitVec{2} \varphi_2(\state_3), \\ 
		h_{4,3}(\state) &= \unitVec{1} \smallint \varphi_3(\state_4) \mathrm{d}\state_4, \quad h_{5,4}(\state) = \unitVec{4} \varphi_4(\state_5), 
	\end{align}
	\end{subequations}
	we have for all $ \state \in \mathbb{R}^5 $
	\begin{subequations}
	\begin{align}
		\unitVec{1} \varphi_1(\state_2) \varphi_2(\state_3) &= [ h_{3,2}, h_{2,1} ](\state) \\
		\unitVec{1} \varphi_3(\state_4) \varphi_4(\state_5) &= [ h_{5,4}, h_{3,2} ](\state).
	\end{align}
	\end{subequations}
	However, we cannot directly use~\eqref{eqGeneralFormula}
	to rewrite a non-admissible vector field of the form $ \unitVec{1} \varphi_2(\state_3) \varphi_4(\state_5) $.
	\oprocend
\end{example}
\colorlet{nodeColor}{green!50!black}
\tikzstyle{nodeStyle} = [circle,draw=white,text=white,fill=nodeColor,thick]
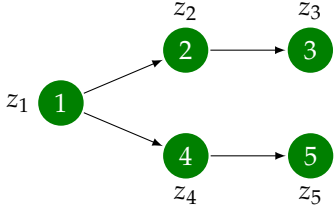
\begin{figure}[t]
	\begin{center}
	\begin{tikzpicture}[>=latex]
		\node[nodeStyle,name=node1] {$1$};
		\node[nodeStyle,name=node2,right=1cm of node1,anchor=west,yshift=20pt] {$2$};
		\node[nodeStyle,name=node3,right=1cm of node2,anchor=west] {$3$};
		\node[nodeStyle,name=node4,right=1cm of node1,anchor=west,yshift=-20pt] {$4$};
		\node[nodeStyle,name=node5,right=1cm of node4,anchor=west] {$5$};
		\draw[->] (node1) -- (node2);
		\draw[->] (node2) -- (node3);
		\draw[->] (node1) -- (node4);
		\draw[->] (node4) -- (node5);
		\node[left=0.2em of node1.west,anchor=east,inner sep=0pt] {$z_1$}; 
		\node[above=0.2em of node2.north,anchor=south,inner sep=0pt] {$z_2$}; 
		\node[above=0.2em of node3.north,anchor=south,inner sep=0pt] {$z_3$}; 
		\node[below=0.2em of node4.south,anchor=north,inner sep=0pt] {$z_4$}; 
		\node[below=0.2em of node5.south,anchor=north,inner sep=0pt] {$z_5$}; 
	\end{tikzpicture}
	\end{center}
	\caption{The graph considered in \cref{exampleGeneralFormula} where the $i$th agent has an associated state $ \state_i $.}\label{figGraphExampleProduct}
\end{figure}
In the next result, we wish to overcome this limitation and 
show how to not only write sums of the vector fields in the form of 
the right-hand side of \eqref{eqGeneralFormula} in terms of 
Lie brackets of admissible vector fields, but also products thereof.
\ifthenelse{\boolean{longVersion}}{
\begin{table*}[t]
\begin{center}
\begin{tabular}{@{}llll@{}}
	\toprule 
	Vector field & & Admissible if... & Rewritable if... \\ \midrule 
	$ \unitVec{i} \tfrac{\partial F_i^{(\ell)}}{\partial x_i}(x_i) \prod\limits_{\substack{j\in\mathcal{J}_F^{(\ell)}\\j\neq i}} F_j^{(\ell)}(x_j) $ & $ \begin{array}{l} \ell = 1,2,\dots,n_F \\ i = 1,2\dots,n \end{array} $ & $ \elLaplacian_{ij} \neq 0 $ for all $ j \in \mathcal{J}_F^{(\ell)} $ & $ \exists p_{ij} $ in $ \mathcal{G} $ for all $ j \in \mathcal{J}_F^{(\ell)} $ \\[1.5em]
	$ \unitVec{i} \bar{a}_{ki} \dualEq_k $ & $ \begin{array}{l} k \in \setEq \\ i = 1,2\dots,n \end{array} $ & $ \elLaplacian_{ik} \neq 0 $ for all $ k: a_{ki} \neq 0 $ & $ \exists p_{ik} $ in $ \mathcal{G} $ for all $ k: a_{ki} \neq 0 $\\[1.5em]
	$ \unitVec{i} \tfrac{\partial c_{k,k}^{(\ell)}}{\partial x_k}(x_k) \prod\limits_{\substack{j\in\mathcal{J}_{c_k}^{(\ell)}\\j\neq k}} c_{k,j}^{(\ell)}(x_j) $ & $ \begin{array}{l} k \in \setIneq, \; \ell = 1,2,\dots,n_{c_k} \\ i =1,2,\dots,n \end{array} $ & $ \elLaplacian_{ij} \neq 0 $ for all $ j \in \mathcal{J}_{c_k}^{(\ell)} $ & $ \exists p_{ij} $ in $ \mathcal{G} $ for all $ j \in \mathcal{J}_{c_k}^{(\ell)} $ \\[1.5em]
	$ \unitVec{n+k} a_{ki} x_i $ & $ \begin{array}{l} k \in \setEq \\ i = 1,2\dots,n \end{array} $ & \Cref{assMatchingConstraints} holds & $ \exists p_{ki} $ in $ \mathcal{G} $ for all $ a_{ki} \neq 0 $ \\[1.5em]
	$ \unitVec{n+\nEq+k} \dualIneq_k \prod\limits_{j\in\mathcal{J}_{c_k}^{(\ell)}} c_{k,j}^{(\ell)}(x_j) $ & $ \begin{array}{l} k \in \setIneq \\ \ell = 1,2,\dots,n_{c_k} \end{array} $ & \Cref{assMatchingConstraints} holds & $ \exists p_{kj} $ in $ \mathcal{G} $ for all $ j \in  \mathcal{J}_{c_k}^{(\ell)} $  \\
	\bottomrule
\end{tabular}
\caption{An overview of all vector fields appearing in~\eqref{eqSPAReformulated}.}\label{tableVectorFields}
\end{center}
\end{table*}
}{}
\begin{proposition}\label{propProduct}
Let $ \eta_k: \mathbb{R} \to \mathbb{R} $, $ \eta_k \in \mathcal{C}^1 $, $ k = 1,2,\dots,m$,
	$m\geq2$,	and define
	\begin{subequations}
	\begin{align}
		\psi_1(\state) &= \unitVec{i} \eta_1(\state_{j_1}) \smallint \eta_0(\state_i) \mathrm{d}\state_i, \\
		\psi_k(\state) &= \unitVec{i} \state_i \eta_k(\state_{j_k}), \quad k = 2,\dots,m-1, \\
		\psi_m(\state) &= \unitVec{i} \eta_m(\state_{j_m}),
	\end{align}\label{eqPropProductDefPsi}
	\end{subequations}
	$ i, j_k \in \mathbb{\dimState} $, $ \psi_k: \mathbb{R}^{\dimState} \to \mathbb{R}^{\dimState} $. Then,
	for any $ j_k \neq i $, $ k~=~1,\dots,m $, we have for all $ \state \in \mathbb{R}^{\dimState} $ 
	\begin{align} 
		\begin{split}
		\bigg[ \Big[ \dots \big[ [ \psi_m, \psi_{m-1} ], \psi_{m-2} \big], \dots \Big], \psi_1 \bigg](\state) \\
		= \unitVec{i} \eta_0(\state_i) \prod\limits_{k=1}^{m} \eta_k(\state_{j_k}).
		\end{split}
		\label{eqPropProduct}
	\end{align}\par\vspace*{-1.5em}\oprocend
\end{proposition}
A proof is given in \ifthenelse{\boolean{longVersion}}{\cref{appProofPropProduct}}{\cite{mic2018extensionsArxiv}}.
The vector fields $ \psi_k $ defined by~\eqref{eqPropProductDefPsi} are 
in general non-admissible; thus, the left-hand side of~\eqref{eqPropProduct}
is not a Lie bracket of admissible vector fields. However, observing that
the vector fields $ \psi_k $ take the same form as the right-hand side
of~\eqref{eqPropFormulaSpecialCase}, we can make use of \cref{propFormulaSpecialCase}
to write the left-hand side of~\eqref{eqPropProduct} as a Lie bracket of
admissible vector fields as long as, for any \mbox{$k\in\lbrace1,\dots,m\rbrace$}, there exists a
path from the $i$th node to the node that is associated to state $ \state_{j_k} $.
We illustrate that by means of \cref{exampleGeneralFormula}.
\setcounter{example}{0}
\begin{example}[continued]
Reconsider \cref{exampleGeneralFormula} and suppose we want to
	rewrite the non-admissible vector field $ \unitVec{1} \varphi_2(\state_3) \varphi_4(\state_5) $.
	Following \cref{propProduct} we let 
	\ifthenelse{\boolean{longVersion}}{
	\begin{align}
		\psi_1(\state) = \unitVec{1} \state_1 \varphi_2(\state_3), \quad \psi_2(\state) = \unitVec{1} \varphi_3(\state_5),
	\end{align}
	}{$\psi_1(\state)~=~\unitVec{1} \state_1 \varphi_2(\state_3)$, $ \psi_2(\state) = \unitVec{1} \varphi_3(\state_5) $,}
	and observe that for all $ \state \in \mathbb{R}^{\dimState} $ we have
	\begin{align}
		[ \psi_2, \psi_1 ](\state) = \unitVec{1} \varphi_2(\state_3) \varphi_4(\state_5) . \label{eqExampleGeneralFormulaContinuedStep1}
	\end{align}
	Further, following \cref{propFormulaSpecialCase} and choosing
	\begin{subequations}
	\begin{align}
		h_{2,1}(\state) &= \unitVec{1} \state_1, \quad && h_{3,2}(\state) = \unitVec{2} \varphi_2(\state_3), \\
		h_{4,3}(\state) &= \unitVec{1},                && h_{5,4}(\state) = \unitVec{4} \varphi_4(\state_5),
	\end{align}
	\end{subequations}
	we have for all $ \state \in \mathbb{R}^\dimState $
	\begin{align}
		\psi_1(\state) = [ h_{3,2}, h_{2,1} ](\state), \quad \psi_2(\state) = [h_{5,4}, h_{4,3}](\state).
	\end{align}
	Using this in~\eqref{eqExampleGeneralFormulaContinuedStep1} we finally managed to
	rewrite $ \unitVec{1} \varphi_2(\state_3) \varphi_4(\state_5) $ in terms 
	of admissible vector fields.
	\oprocend
\end{example}
\begin{remark}
Instead of realizing the multiplication by means of Lie brackets another way
	is to augment the agent state by estimates of the respective state of
	the other agent. More precisely, in~\cref{exampleGeneralFormula}, we augment
	the state of agent~1 by $ \xi_3 $ and $ \xi_5 $ that are estimates of $ \state_3 $
	and $ \state_5 $, respectively, and let
	\begin{align}
		\begin{bmatrix}
		\dot{\state}_1 \\ \dot{\xi}_3 \\ \dot{\xi}_5
		\end{bmatrix}
		&=
		\begin{bmatrix}
		w(\state_1,\state_2,\state_4,\xi_3,\xi_5) \\
		- \mu \xi_3 \\
		- \mu \xi_5
		\end{bmatrix}
		+ \unitVec{2} \state_5 + \unitVec{3} \state_5,
	\end{align}
	where $ w: \mathbb{R}^5 \to \mathbb{R} $ and $ \mu > 0 $ is sufficiently large, hence $ \xi_3(t) \approx \state_3(t) $,
	$ \xi_5(t) \approx \state_5(t) $. The resulting non-admissible vector fields in the complete augmented system
	can then be written in terms of Lie brackets of admissible vector fields
	using~\cref{propFormulaSpecialCase}. However, in the application at
	hand this alters the saddle-point dynamics~\eqref{eqSPA} which necessitates
	a stability analysis of the augmented system.
\end{remark}
Hence, under suitable assumptions on the communication graph,
this allows us to write vector fields whose components are 
sums of products of arbitrary functions in terms of Lie brackets
of admissible vector fields. This observation gives rise to the 
next Lemma.
\begin{lemma}\label{lemmaAnalyticFunctions}
Consider a strongly connected graph $ \mathcal{G} $ of $n$ nodes
	and let $ \varphi: \mathbb{R}^{\dimState} \to \mathbb{R} $ be an
	analytic function.
	Then any vector field $ \psi(\state) = \unitVec{i} \varphi(\state) $, $ \unitVec{i} \in \mathbb{R}^{\dimState} $, $ i = 1,2,\dots,\dimState $,
	can be written as a possibly infinite sum of Lie brackets
	of admissible vector fields.
	\oprocend
\end{lemma}
\ifthenelse{\boolean{longVersion}}{
\begin{proof}
Since $ \varphi $ is analytic, by a {series} expansion it can be
	written as a possibly infinite sum of monomials of the components $ \state_j $
	of $ \state $, $ j = 1,2,\dots,\dimState$. Using \cref{propProduct},
	all these monomials can be written in terms of Lie brackets
	of vector fields of the form~\eqref{eqPropProductDefPsi}.
	By strong connectivity of $\mathcal{G}$, all these vector fields
	can be written in terms of Lie brackets of admissible vector fields,
	thus concluding the proof.
\end{proof}
}{A proof is given in~\cite{mic2018extensionsArxiv}.}
While the result might be more of a theoretical nature for the 
application at hand, it nevertheless shows that the proposed
approach in principle applies to a large class of problems.

\subsection{Distributed optimization via Lie brackets} 
\def\deltaAngle{20}
\begin{figure*}[t]
	\begin{center}
	\begin{floatrow}
	\ffigbox[\FBwidth]{
		\begin{tikzpicture}[>=latex]
			\node[nodeStyle,draw,name=node1,anchor=center] at(90+36:1.5) {$1$};
			\node[nodeStyle,draw,name=node2,anchor=center] at(90-36:1.5) {$2$};
			\node[nodeStyle,draw,name=node3,anchor=center] at(90-36-72:1.5) {$3$};
			\node[nodeStyle,draw,name=node4,anchor=center] at(90-36-72-72:1.5) {$4$};
			\node[nodeStyle,draw,name=node5,anchor=center] at(90-36-72-72-72:1.5) {$5$};
			\draw[->] (node1) -- (node2);
			\draw[->] (node2) -- (node3);
			\draw[->] (node2) -- (node5);
			\draw[->] (node3) -- (node4);
			\draw[->] (node4) -- (node5);
			\draw[->] (node5) -- (node1);
			\node[below=0.7cm of node4,minimum width=4cm] {};
		\end{tikzpicture}
	}{
		\caption{The communication graph for the example considered in~\cref{secExample}.}\label{figGraphExample}
	}
	\capbtabbox{
		\begin{tabular}{@{}llll@{}}
			\toprule
			Vector field & Corresponding path & Lie bracket representation   & \\ \midrule
			$ \unitVec{5} \dualEq_2 $ & $ p_{52} = \pathLeft v_5 \pathSep v_1 \pathSep v_2 \pathRight $ & $ \big[ \unitVec{j} \state_6, \unitVec{5} \state_j \big] $ & $ j \in \mathcal{I}(1) $ \\[0.3em]
			$ \unitVec{2} x_2 \dualIneq_1 $ & $ p_{21} = \pathLeft v_2 \pathSep v_5 \pathSep v_1 \pathRight $ & $ \big[ \unitVec{j} \state_7, \unitVec{2} \state_2 \state_j \big] $ & $ j \in \mathcal{I}(5) $ \\[0.3em]
			$ \unitVec{1} \dualIneq_3 $ & $ p_{13} = \pathLeft v_1 \pathSep v_2 \pathSep v_3 \pathRight $ & $ \big[\unitVec{j} \dualIneq_3, \unitVec{1} \state_{j} \big] $ & $ j \in \mathcal{I}(2) $ \\[0.3em]
			$ \unitVec{4} x_1 \dualIneq_4 $ & $ p_{41} = \pathLeft v_4 \pathSep v_5 \pathSep v_1 \pathRight $ & $ \big[ \unitVec{j} x_1, \unitVec{4} \dualIneq_4 \state_j \big] $ & $ j \in \mathcal{I}(5) $ \\[0.3em]
			$ \unitVec{1} x_4 \dualIneq_4 $ & $ p_{14} = \pathLeft v_1 \pathSep v_2 \pathSep v_3 \pathSep v_4 \pathRight $ & $ \big[ \unitVec{j_2} \dualIneq_4 x_4, [ \unitVec{j_1} \state_{j_2}, \unitVec{1} \state_{j_1} ] \big] $ & $ j_1 \in \mathcal{I}(2) $, $ j_2 \in \mathcal{I}(3) $ \\[0.3em]
			$ \unitVec{1} x_1 \dualIneq_4 $ & $ p_{14} = \pathLeft v_1 \pathSep v_2 \pathSep v_3 \pathSep v_4 \pathRight $ & $ \big[ \unitVec{j_2} \dualIneq_4, [\unitVec{j_1} \state_{j_2}, \unitVec{1} x_1 \state_{j_1} ] \big] $ & $ j_1 \in \mathcal{I}(2) $, $ j_2 \in \mathcal{I}(3) $ \\[0.3em]
			$ \unitVec{8} \dualIneq_3 x_1 $ & $ p_{31} = \pathLeft v_3 \pathSep v_4 \pathSep v_5 \pathSep v_1 \pathRight $ & $ \big[ \unitVec{j_2} x_1, [ \unitVec{j_1} \state_{j_2}, \unitVec{8} \dualIneq_2 \state_{j_1} ] \big] $ & $ j_1 \in \mathcal{I}(4) $, $ j_2 \in \mathcal{I}(5) $ \\[0.3em]
			$ \unitVec{9} \dualIneq_4 x_4 x_1 $ & $ p_{41} = \pathLeft v_4 \pathSep v_5 \pathSep v_1 \pathRight $ & $ \big[ \unitVec{j}x_1, \unitVec{9} \state_j x_4 \dualIneq_4 ] \big] $ & $ j \in \mathcal{I}(5) $ \\[0.3em]
			$ \unitVec{9} \dualIneq_4 x_1^2 $ & $ p_{41} = \pathLeft v_4 \pathSep v_5 \pathSep v_1 \pathRight $ & $ \big[ \unitVec{j} x_1^2, \unitVec{9} \dualIneq_4 \state_j \big] $ & $ j \in \mathcal{I}(5) $ \\
			\bottomrule
		\end{tabular}
	}{
	\caption{An overview of the non-admissible vector fields in~\eqref{eqSPAExample} and their Lie bracket representations.
	Here, the index sets are $ \mathcal{I}(1) = \lbrace 1,7 \rbrace $, $ \mathcal{I}(2) = \lbrace 2, 6 \rbrace $, $ \mathcal{I}(3) = \lbrace 3,8 \rbrace $,
	$ \mathcal{I}(4)=\lbrace 4,9 \rbrace $, $ \mathcal{I}(5) = \lbrace 5 \rbrace $.
	}\label{tableExampleNonAdmissibleVectorFields}
	}
	\end{floatrow}
	\end{center}
	\ifthenelse{\boolean{longVersion}}{}{\vspace*{-1em}}
\end{figure*}
In the sequel we apply the results from the last section
to the problem at hand and rewrite the saddle-point dynamics~\eqref{eqSPA}
by means of Lie brackets of admissible vector fields. 
For the sake of a simpler notation we assume in the following that
$ {\nEq}_i = 1 $ for all $ i \in \setEq $ and $ {\nIneq}_i = 1 $
for all $ i \in \setIneq $. We further assume that each agent
has an associated equality and inequality constraint,
i.e., $ \setEq ~=~\lbrace 1,2,\dots,n \rbrace $,
$ \setIneq~=~\lbrace 1,2,\dots,n \rbrace $. This can always be
achieved by augmenting the optimization problem~\eqref{eqOptimizationProblem}
by constraints that do not alter the feasible set. We emphasize that this
is not necessary for the methodology to apply as we will illustrate
in the example in~\cref{secExample}. We can then write the saddle-point 
dynamics~\eqref{eqSPA} equivalently as
\begin{subequations}
\begin{flalign}
	\dot{x} &= - \sum\limits_{i=1}^{n} \unitVec{i} \big( \tfrac{\partial F}{\partial x_i}(x) + \sum\limits_{k\in\setEq} \hspace*{-0.3em} \bar{a}_{ik} \dualEq_k + \sum\limits_{k\in\setIneq} \hspace*{-0.5em} \tfrac{\partial c_k}{\partial x_k}(x) \dualIneq_k \big) \hspace*{-2em} & \label{eqSPAReformulatedA} \\
	\dot{\dualEq} &= \sum\limits_{k \in \setEq}\unitVec{k} \big( \sum\limits_{i=1}^{n}  \bar{a}_{ki} x_i + \bar{b}_k \big) \label{eqSPAReformulatedB} \\
	\dot{\dualIneq} &= \sum\limits_{k \in \setIneq} \unitVec{k} \dualIneq_{k} c_k(x), \label{eqSPAReformulatedC}
\end{flalign}\label{eqSPAReformulated}
\end{subequations}
where $ a_i(x) = \bar{a}_i x + \bar{b}_i $, $ \bar{a}_{i} =  [ \bar{a}_{i1}, \bar{a}_{i2}, \dots, \bar{a}_{in} ] \in \mathbb{R}^{1 \times n} $, 
$ \bar{b}_i \in \mathbb{R} $. 
{Motivated by our previous discussions, we assume
in the following} that the objective function as well
as the inequality constraints are sums of products of separable
functions, i.e.,
\begin{align}
	F(x) &= \sum\limits_{\ell=1}^{n_F}  \prod\limits_{j \in \mathcal{J}_F^{(\ell)} } F_j^{(\ell)}(x_j), \label{eqStructureObjectiveFunction}  \\
	c_k(x) &= \sum\limits_{\ell=1}^{n_{c_k}}  \prod\limits_{j \in \mathcal{J}_{c_k}^{(\ell)} } c_{k,j}^{(\ell)}(x_j), \label{eqStructureInequalityConstraint}
\end{align}
where the $ F_j^{(\ell)}: \mathbb{R} \to \mathbb{R} $, $ j\in\mathcal{J}_F^{(\ell)} {\subseteq \lbrace 1,2,\dots,n\rbrace }$,
\mbox{$\ell=1,2,\dots,n_F$}, $n_F~>~0$, are strictly convex functions and the 
$ c_{k,j}^{(\ell)}: \mathbb{R} \to \mathbb{R} $, $ k\in\setIneq $,
$j\in\mathcal{J}_{c_k}^{(\ell)} {\subseteq \lbrace 1,2,\dots,n\rbrace }$, $\ell=1,2,\dots,n_{c_i}$, $n_{c_k}>0$, are convex.
Observe that, if $ n_F $ and $ n_{c_k} $ are infinite and $ F_j^{(\ell)} $,
$ c_{k,j}^{(\ell)} $ are monomials, this includes all analytic functions $ F $, $ c_k $.
If $ \mathcal{J}_F^{(\ell)} = \lbrace \ell \rbrace $, $ \mathcal{J}_{c_k}^{(\ell)} = \lbrace \ell \rbrace $,
we obtain the particularly important
special case that both the objective function and the constraints 
are a sum of separable functions; hence also
the case of linear constraints considered in~\cite{ebenbauer2017directed},
\cite{mic2017extremum} is covered here.
\ifthenelse{\boolean{longVersion}}{Under this assumption the vector fields appearing in~\eqref{eqSPAReformulated}
are summed up in~\cref{tableVectorFields}.
Depending on the communication graph as well as the structure
of the constraints and the objective function, these
vector fields can either be admissible or not.}
{Depending on the communication graph as well as the structure
of the constraints and the objective function, the
vector fields in~\eqref{eqSPAReformulated} can either be admissible or not.}
In particular,
the vector fields in~\eqref{eqSPAReformulatedB}, \eqref{eqSPAReformulatedC}
are admissible if the constraints are compatible with the communication
topology defined by the graph $ \mathcal{G} $, i.e., if the following
assumption holds: 
\begin{assumption}\label{assMatchingConstraints}
For all $ i,j = 1,2,\dots,n $ with $ \elLaplacian_{ij} = 0 $ we have 
	$ \tfrac{\partial a_i}{\partial x_j}(x) \equiv 0 $
	as well as $ \tfrac{\partial c_i}{\partial x_j}(x) \equiv 0 $.
	\oprocend
\end{assumption}
We point out that all non-admissible vector fields in~\eqref{eqSPAReformulated}
can be written in terms of Lie brackets of admissible vector fields
under appropriate assumptions on the communication graph, see also \ifthenelse{\boolean{longVersion}}{the last column 
of \cref{tableVectorFields}}{\cite[Table 1]{mic2018extensionsArxiv}}. Specifically, if the graph is strongly
connected, then all non-admissible vector fields can be rewritten
independent of the objective function as well as the constraints,
given that they admit the structure~\eqref{eqStructureObjectiveFunction},
\eqref{eqStructureInequalityConstraint}. In most cases, however, 
much less restrictive requirements on the communication graph
are sufficient. We do not explicitly discuss how to rewrite the 
non-admissible vector fields using \cref{propFormulaSpecialCase} 
and \cref{propProduct} in general but illustrate this by means of an example in 
\cref{secExample}.

\section{EXAMPLE}\label{secExample}
In this section we illustrate the previous results by
means of an example. Consider the following optimization 
problem
\begin{align}
	\begin{split}
	\min\limits_{x\in\mathbb{R}^5}        \quad & F(x) = \sum\limits_{i=1}^{5} F_i(x_i) \\
	\text{s.t.} \quad & a_2(x) = 2 x_2 - x_5 = 0 \\
	                  & c_1(x) = x_1^2 + x_2^2 - 4 \leq 0 \\
	                  & c_3(x) = x_1 + x_3 - 2 \leq 0 \\
	                  & c_4(x) = x_4^2 - x_4x_1 + x_1^2 - {100} \leq 0,
	\end{split}\label{eqExampleOptProb}
\end{align}
where $ F_i(x_i) = ( x_i - i )^2 $, $ x = [ x_1, x_2, x_3, x_4, x_5 ]^\top \in \mathbb{R}^5 $.
We assume that the communication topology is described by
the graph in~\cref{figGraphExample}. Observe that the constraints
$ c_3 $ and $ c_4 $ are not compatible with the graph topology, hence
\cref{assMatchingConstraints} is not fulfilled.
The corresponding saddle-point dynamics~\eqref{eqSPA} are then given by
\begin{subequations}
\begin{flalign}
	\dot{x} =\;&-\nabla F(x) - ( 2 \unitVec{2} - \unitVec{5} ) \dualEq_2 - 2( x_1\unitVec{1} + x_2\unitVec{2} ) \dualIneq_1 \hspace*{-2em} & \label{eqSPAExampleA}  \\
	&- ( 2x_4 - x_1 ) \unitVec{4} \dualIneq_4 - (2x_1 - x_4) \unitVec{1} \dualIneq_4 - (\unitVec{1} + \unitVec{3} ) \dualIneq_3 \nonumber \\
	\dot{\dualEq}_2 =\;&2x_2 - x_5 \label{eqSPAExampleB} \\
	\dot{\dualIneq} =\;&\dualIneq_1 \unitVec{1} ( x_1^2 + x_2^2 - 4 ) + \dualIneq_3 \unitVec{2} ( x_1 + x_3 - 2 ) \nonumber  \\
	&+ \dualIneq_4 \unitVec{3} ( x_4^2 - x_4 x_1 + x_1^2 - {100} ) , \label{eqSPAExampleC} 
\end{flalign}\label{eqSPAExample}
\end{subequations}
where $ x \in \mathbb{R}^5 $, $ \dualEq_2 \in \mathbb{R} $, and 
$ \dualIneq=[ \dualIneq_1, \dualIneq_3, \dualIneq_4 ]^\top \in \mathbb{R}^3 $.
We next rewrite all non-admissible vector fields in~\eqref{eqSPAExample}
using \cref{propFormulaSpecialCase}. 
We will thereby follow the choice~\eqref{eqVectorFieldsSimpleChoice} to make
sure that~\eqref{eqPropConstraintFunctions} holds. We do not discuss
how to rewrite each non-admissible vector field in detail, but limit
ourselves to the vector field $ \unitVec{1} x_4 \dualIneq_4 $ 
from~\eqref{eqSPAExampleA}. Comparing with~\eqref{eqPropFormulaSpecialCase},
we have $ j_1 = 1 $, $ j_r = 4 $, and $ f_{j_1}^{(1)}(\state_{\mathcal{I}(j_1)}) = 1 $,
$ f_{j_{r-1}}^{(2)}( \state_{\mathcal{I}(j_r)} ) = x_4 \dualIneq_4 $.
Following~\cref{propFormulaSpecialCase},
we require a path from node~$1$ to node~$4$ that is here given by
$ p_{14} = \pathLeft v_1 \pathSep v_2 \pathSep v_3 \pathSep v_4 \pathRight $; thus $ r = 4 $.
We then obtain
\ifthenelse{\boolean{longVersion}}{
\begin{align}
	\unitVec{1} x_4 \dualIneq_4 = \big[ \unitVec{j_2} \dualIneq_4 x_4, [ \unitVec{j_1} \state_{j_2}, \unitVec{1} \state_{j_1} ] \big],
\end{align}
}{$ \unitVec{1} x_4 \dualIneq_4 = \big[ \unitVec{j_2} \dualIneq_4 x_4, [ \unitVec{j_1} \state_{j_2}, \unitVec{1} \state_{j_1} ] \big] $,}
where $ j_1 \in \mathcal{I}(2) = \lbrace 2,6 \rbrace $, $ j_2 \in \mathcal{I}(3) = \lbrace 3, 8 \rbrace $,
and the right-hand side is a Lie bracket of admissible vector fields.
All other non-admissible vector fields
in~\eqref{eqSPAExample} can be treated similarly and we sum up the resulting
Lie brackets in~\cref{tableExampleNonAdmissibleVectorFields}. 
For the simulation we let each $ j, j_1, j_2 $ be the largest
index in its respective index set. By that choice, we inject less
perturbation in the primal and more in the dual variables New{
which is also visible in the simulation results depicted in~\cref{figExampleSimRes}.}
{As to be seen, the distributed algorithm approximates the trajectories
of the non-distributed saddle-point dynamics~\eqref{eqSPAExample}
and converges to a neighborhood of the optimizer~$ x^\star = [ 0, 2, 2, 4, 4 ]^\top $.}
New{We also included simulation results with additional low-pass filters
in the distributed $x$-, $\dualEq$- and $\dualIneq$-dynamics. 
While the effect on the primal variables is small since we already reduced
the oscillations by our design choice, the dual variables show significantly
less oscillations and better approximate the non-distributed trajectories.
A rigorous stability analysis of the augmented distributed dynamics and 
a performance-oriented design of the filters is up to future work.}
\ifthenelse{\boolean{longVersion}}{
\begin{figure}[t]
	\centering
	\includegraphics[width=\textwidth]{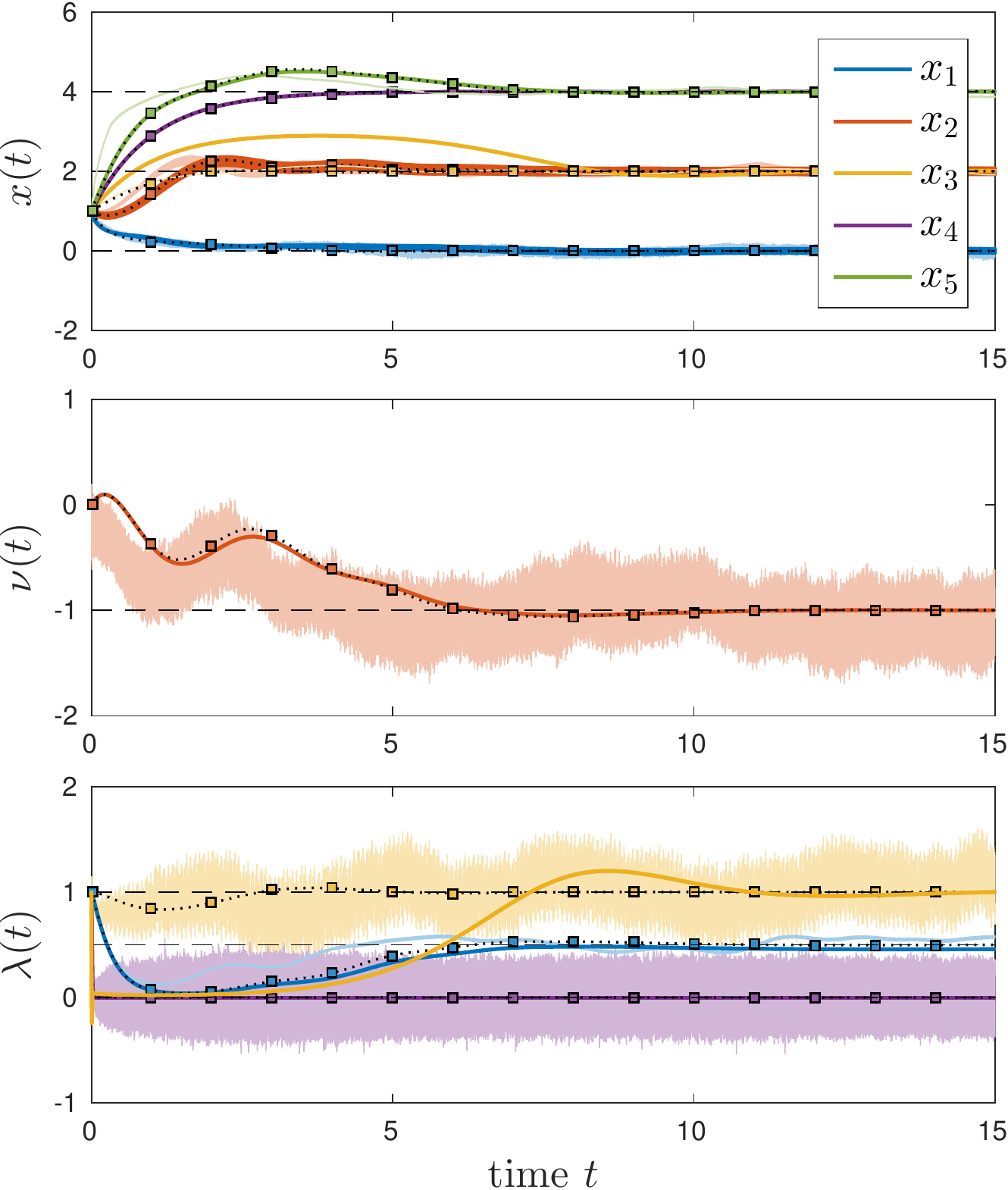}
	\caption{Simulation results of the example from~\cref{secExample}. The upper plot shows the primal variable~$x(t)$,
	the lower left one the dual variable~$\dualEq_2(t)$ and the lower right one the dual variable~$\dualIneq(t)$.
	New{In each plot, the dotted black lines marked with squares depict the trajectories of the
	non-distributed saddle-point dynamics~\eqref{eqSPAExample} while the thinner oscillating ones
	depict the trajectories of the distributed approximation. The corresponding thick lines
	depict the distributed approximation with additional low-pass filters. The dashed black lines
	indicate the desired equilbrium of~\eqref{eqSPAExample}.}
	}\label{figExampleSimRes}
\end{figure}}{
\begin{figure}[t]
	\centering
	\includegraphics[width=\textwidth]{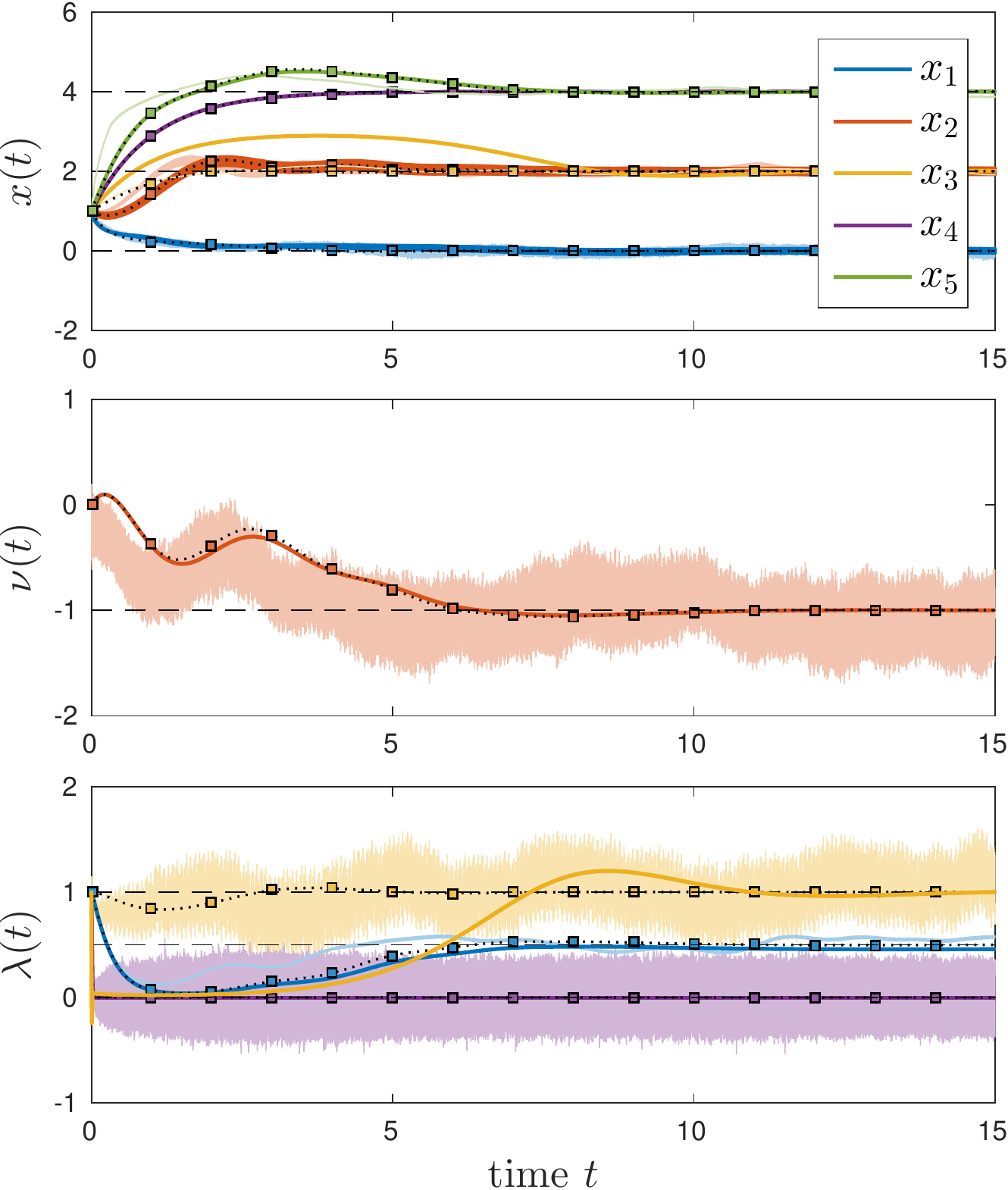}
	\caption{Simulation results of the example from~\cref{secExample}. 
	The upper plot shows the primal variable~$x(t)$,
	the middle one the dual variable~$\dualEq_2(t)$ and the lower one the dual variable~$\dualIneq(t)$.
	New{In each plot, the dotted black lines marked with squares depict the trajectories of the
	non-distributed saddle-point dynamics~\eqref{eqSPAExample} while the thinner oscillating ones
	depict the trajectories of the distributed approximation. The corresponding thick lines
	depict the distributed approximation with additional low-pass filters. The dashed black lines
	indicate the desired equilbrium of~\eqref{eqSPAExample}.}
	}\label{figExampleSimRes}
	\vspace*{-1em}
\end{figure}}

\section{CONCLUSIONS AND OUTLOOK}\label{secConclusions}
We considered a convex optimization problem and showed
how distributed optimization algorithms can be designed
for a quite general class of problems with little structural
requirements under mild assumptions on the communication network.
We therefore extended the Lie bracket approximation approach to
distributed optimization proposed in~\cite{ebenbauer2017directed},
\cite{mic2017extremum} and discussed which kind of vector fields
can in principle be written in terms of Lie brackets of admissible
vector fields. We did not discuss the construction of approximating
inputs {but postpone this to~\cite{mic2017opt} where we will
present a modified version of the general algorithm from~\cite{liu1997approximation}
that exploits the structure of the problem at hand.}

\bibliographystyle{IEEEtran}
\bibliography{IEEEabrv,subfiles/bibfile}

\section{APPENDIX}\label{secAppendix}
\subsection{Proof of \Cref{lemmaGeneralFormula}}\label{appProofLemmaGeneralFormula}
\begin{proof}
We prove this result by induction. First, for $ r = 1 $, \eqref{eqGeneralFormula} is
	trivially true by definition~\eqref{eqDefVectorFields}. Suppose now that
	the claim holds for all $ r \leq \bar{r} $, $ \bar{r} > 1 $ and consider $ r = \bar{r} + 1 $.
	Define 
	\begin{align}
		\begin{split}
		&\tilde{f}_{j_{\bar{r}},j_1}(\state) =\\ 
		& f_{j_{\bar{r}-1}}( \state_{\mathcal{I}(j_{\bar{r}-1})}, \state_{\mathcal{I}(j_{\bar{r}})} ) \prod\limits_{k=1}^{\bar{r}-2} \tfrac{\partial f_{j_k}}{\partial \state_{j_{k+1}}}(\state_{\mathcal{I}(j_k)},\state_{\mathcal{I}(j_{k+1})}).
		\end{split}
		\label{eqDefFTildeProof}
	\end{align}
	By the induction hypothesis we then have
	\begin{align}
		 &~ \bigg[ h_{j_{\bar{r}+1},j_{\bar{r}}}, \Big[ h_{j_{\bar{r}},j_{\bar{r}-1}}, \big[ \dots, [ h_{j_{3},j_{2}}, h_{j_{2},j_{1}} ] \dots \big] \Big] \bigg](\state) \nonumber \\
		=&~ [ h_{j_{\bar{r}+1},j_{\bar{r}}}, \unitVec{j_1} \tilde{f}_{j_{\bar{r}},j_1}  ](\state) \nonumber \\
		=&~ \unitVec{j_1} \big( \sum\limits_{i=1}^{\dimState} \tfrac{\partial \tilde{f}_{j_{\bar{r}},j_1}}{\partial \state_{i}}(\state) \big) \unitVec{j_{\bar{r}}} f_{j_{\bar{r}}} (\state_{\mathcal{I}(j_{\bar{r}})}, \state_{\mathcal{I}(j_{\bar{r}+1})} ) \nonumber \\	
		-&~ \unitVec{j_{\bar{r}}} \big( \sum\limits_{ \substack{ i \in \mathcal{I}(j_{\bar{r}}) \\ \cup \mathcal{I}(j_{\bar{r}+1}) } } \unitVec{i} \tfrac{\partial f_{j_{\bar{r}}} }{\partial \state_i}(\state_{\mathcal{I}(j_{\bar{r}})}, \state_{\mathcal{I}(j_{\bar{r}+1})}) \big) \unitVec{j_1} \tilde{f}_{j_{\bar{r}},j_1}(\state) \nonumber \\
		=&~ \unitVec{j_1} \tfrac{ \partial \tilde{f}_{j_{\bar{r}},j_1} }{ \partial \state_{j_{\bar{r}}} } (\state) f_{j_{\bar{r}}} (\state_{\mathcal{I}(j_{\bar{r}})}, \state_{\mathcal{I}(j_{\bar{r}+1})} ),
	\end{align}
	where $ \unitVec{j_i} \in \mathbb{R}^{n + \nEq + \nIneq} $ is the $ j_i$th unit vector.
	Here, we used that $ p_{i_1 i_{\bar{r}+1}} $ is a simple path; hence, since
	$ \mathcal{I}(i_{k_1}) $ and $ \mathcal{I}(i_{k_2}) $ are disjunct for any
	$ k_1 \neq k_2 $, also $ j_{k_1} \neq j_{k_2} $ for any $ k_1, k_2 = 1,2,\dots, \bar{r} + 1 $,
	$ k_1 \neq k_2 $. Using definition~\eqref{eqDefFTildeProof}, we obtain~\eqref{eqGeneralFormula}.
	Finally, since $ j_k \in \mathcal{I}(i_k) $ and $ g_{i_k i_{k+1}} \neq 0 $ for
	all $ k = 1,2,\dots,r-1 $, all $ h_{j_{k+1}, j_{k}} $ are admissible; hence
	the left-hand side of~\eqref{eqGeneralFormula} is a Lie bracket of admissible vector fields.
\end{proof}

\subsection{Proof of \Cref{propFormulaSpecialCase}}\label{appProofPropFormulaSpecialCase}
\begin{proof}
Observe first that
	\begin{align}
		  & f_{j_{r-1}}( \state_{\mathcal{I}(j_{r-1})}, \state_{\mathcal{I}(j_r)} ) \prod\limits_{k=1}^{r-2} \tfrac{\partial f_{j_k}}{\partial \state_{j_{k+1}}}(\state_{\mathcal{I}(j_k)},\state_{\mathcal{I}(j_{k+1})}) \nonumber \\
		=~& f_{j_{r-1}}^{(2)}( \state_{\mathcal{I}(j_{r})}) f_{j_{r-1}}^{(1)}( \state_{\mathcal{I}(j_{r-1})}) \nonumber \\
		& \phantom{f_{j_{r-1}}^{(2)}( \state_{\mathcal{I}(j_{r})})} \times \prod\limits_{k=1}^{r-2} f_{j_k}^{(1)}(\state_{\mathcal{I}(j_k)}) \prod\limits_{k=2}^{r-1} \tfrac{ \partial f_{j_{k-1}}^{(2)} }{ \partial \state_{j_{k}} }(\state_{\mathcal{I}(j_k)}) \nonumber \\
		=~& f_{j_{r-1}}^{(2)}( \state_{\mathcal{I}(j_{r})}) f_{j_1}^{(1)}(\state_{\mathcal{I}(j_1)}) \prod\limits_{k=1}^{r-1} f_{j_k}^{(1)}(\state_{\mathcal{I}(j_k)}) \tfrac{ \partial f_{j_{k-1}}^{(2)} }{ \partial \state_{j_{k}} }(\state_{\mathcal{I}(j_k)}). \nonumber
	\end{align}
	Then, using~\eqref{eqPropConstraintFunctions} and applying \cref{lemmaGeneralFormula},
	the result immediately follows.
\end{proof}

\subsection{Proof of \Cref{propProduct}}\label{appProofPropProduct}
\begin{proof}
We first show by induction that
	\begin{align}
		\bigg[ \Big[ \dots \big[ [ \psi_m, \psi_{m-1} ], \psi_{m-2} \big], \dots \Big], \psi_2 \bigg](\state) 
		= \unitVec{i} \prod\limits_{k=2}^{m} \eta_k(\state_{j_k}). \label{eqProofPropProductIntermediate}
	\end{align}
	For $m=2$ it is clear that \eqref{eqProofPropProductIntermediate} holds.
	Suppose now that~\eqref{eqProofPropProductIntermediate} holds for all $m \leq \bar{m} $,
	$\bar{m} \geq 2 $. For $ m = \bar{m} + 1 $ we then have with a slight abuse
	of notation
	\begin{align}
		&\bigg[ \Big[ \dots \big[ [ \psi_{\bar{m}+1}, \psi_{\bar{m}} ], \psi_{\bar{m}-1} \big], \dots \Big], \psi_2 \bigg](\state) \nonumber \\
		=~& \big[  \unitVec{i} \prod\limits_{k=2}^{\bar{m}+1} \eta_k(\state_{j_k}), \unitVec{i} \state_i \eta_2(\state_{j_2}) \big] \nonumber \\
		=~& \unitVec{i} \big( \unitVec{i}^\top \eta_2(\state_{j_2}) + \state_i \unitVec{j_2} \tfrac{\partial \eta_2}{\partial \state_{j_2}} (\state_{j_2}) \big) \unitVec{i} \prod\limits_{k=2}^{\bar{m}+1} \eta_k(\state_{j_k}) \nonumber \\
		&+ \unitVec{i} \big( \sum\limits_{k=2}^{\bar{m}+1} \unitVec{j_k} \tfrac{\partial \eta_k}{\partial \state_{j_k}}(\state_{j_k}) \prod\limits_{ \substack{ \ell=2 \\ \ell \neq k} } \eta_{\ell}(\state_{j_{\ell}}) \big) \unitVec{i} \state_i \eta_2(\state_{j_2}) \nonumber \\
		=~& \unitVec{i} \prod\limits_{k=2}^{\bar{m}+1} \eta_k(\state_{j_k}),
	\end{align}
	where we used that $ i \neq j_k $ for any $ k = 1,\dots,m $
	in the last step. This proves~\eqref{eqProofPropProductIntermediate}.
	Finally, using~\eqref{eqProofPropProductIntermediate}, we obtain
	\begin{align}
		&\bigg[ \Big[ \dots \big[ [ \psi_m, \psi_{m-1} ], \psi_{m-2} \big], \dots \Big], \psi_1 \bigg](\state) \nonumber \\
		=~& \big[ \unitVec{i} \prod\limits_{k=2}^{m} \eta_k(\state_{j_k}), \unitVec{i} \eta_1(\state_{j_1}) \smallint \eta_0(\state_i) \mathrm{d}\state_i \big] \nonumber \\
		=~& \unitVec{i} \eta_0(\state_i) \prod\limits_{k=1}^{m} \eta_k(\state_{j_k}).
	\end{align}
	This finishes the proof.

\end{proof}

\end{document}